%Paper: dg-ga/9406003
%From: Mike Wolf <mwolf@msri.org>
%Date: Wed, 29 Jun 1994 09:44:26 -0700

\documentstyle{amsppt}
%\pagewidth{6in}
%\pageheight{8in}
\NoRunningHeads
\NoBlackBoxes
\magnification = 1200

\define\a{\alpha}
\predefine\barunder{\b}
\redefine\b{\beta}
\define\bpage{\bigpagebreak}

\predefine\dotunder{\d}
\redefine\d{\delta}
\define\e{\epsilon}

\define\g{\gamma}
\define\G{\Gamma}
\define\k{\kappa}
\define\lb\{{\left\{}

\define\lm{\limits}

\define\oper{\operatorname}

\define\ov{\overline}

\define\ox{\otimes}
\define\p{\partial}
\define\rb\}{\right\}}

\define\sub{\subheading}

\define\z{\zeta}
\define\({\left(}
\define\){\right)}
\define\[{\left[}
\define\]{\right]}
\define\<{\left<}
\define\>{\right>}
\def\slantline#1#2#3#4#5{\hbox to 0pt{
}}

%SCRIPT LETTERS
%%%%%%%%%%

\def\SF{{\Cal F}}

\def\SR{{\Cal R}}
\def\SS{{\Cal S}}

\def\SU{{\Cal U}}

%%%%%%%
%%%BOLDFACE

\def\BR{{\bold R}}

%%%%%BLACKBOARD BOLD LETTERS

\define\vecs#1#2{#1_1,\dots,#1_#2}
\define\Ree{\oper{Re}}
\def\const{\oper{const}}
\define\im{\oper{Im}}

\topmatter
\title On the Existence of Jenkins-Strebel Differentials Using
Harmonic Maps from Surfaces to Graphs
\endtitle
\author Michael Wolf*\endauthor
\affil
Department of Mathematics\\
Rice University\\
Houston, TX 77251\endaffil
\thanks *Partially supported by NSF grants DMS9300001 and
DMS9022140 (MSRI); Alfred P\. Sloan Research
Fellow\endthanks
\endtopmatter

\document

\centerline{Address until 8/1/94: Mathematical Sciences Research
Institute}

\hskip147pt 1000 Centennial Drive

\hskip147pt  Berkeley, CA 94720

\vskip .7cm
\flushpar AMS Classification Code (1991): 30C70, 30C75, 58E20
\vskip.7cm
\sub{Abstract} We give a new proof of the existence (\cite{HM}, \cite{Ren})
of a Jenkins-Strebel
differential $\Phi$ on a Riemann surface $\SR$ with
prescribed heights of cylinders
by considering the harmonic map from $\SR$ to the leaf space of the
vertical foliation of $\Phi$, thought of as a Riemannian graph.  The novelty
of the argument is that it is essentially Riemannian as well as elementary;
moreover, the harmonic maps existence theory on which it relies is classical,
due mostly to Morrey (\cite{Mo}).
\vskip1cm

\sub{\S1 Introduction} In a series of pathbreaking papers
(\cite{Je}, \cite{Str1}, \cite{Str2}) in the 1950's and 1960's,
Jenkins and Strebel proved the existence of
holomorphic quadratic differentials on a Riemann surface $\SR$, the
complement of whose critical trajectories were Euclidean cylinders
foliated by closed curves; the metrics associated to
these differentials uniquely solved natural extremal
length problems (\cite{Je}), or the free homotopy classes of the
core curves of the cylinders and the ratio of the moduli of the
cylinders could be uniquely specified (\cite{Str1}).
Later, Renelt (\cite{Ren}) and Hubbard and Masur (\cite{HM})
showed that the homotopy classes of these
core curves and the heights of the cylinders
could also be specified. The goal of this note is to
provide another proof of the existence of these differentials with
prescribed heights on $\SR$ based on energy-minimizing maps from
$\SR$ to graphs.

Our method is somewhat unusual in the subject of quadratic
differentials on Riemann surfaces in that our techniques are
essentially Riemannian, with the basic existence theory due mostly
to Morrey \cite{Mo} in 1948; the conformal type of the Riemann
surface $\SR$ is involved because of the conformal invariance of the
total energy of the map. Our holomorphic quadratic differential
on $\SR$ is
the Hopf differential, known classically (and emphasized recently
by Schoen \cite{S}) to result from a stationary point $I:\SR\to N$
of a conformally invariant functional (here, it is important that
$I$ be stationary with respect to reparametrizations of the domain
$\SR$). A novel feature of our argument is its use of maps with domain
$\SR$, rather than, say, maps of cylinders into a range Riemann
surface $\SR$.

Here, roughly, is the proof. Draw the foliation with $k$ prescribed
cylinders on $\SR$, and suppose for now that the leaf space is a
boundaryless graph $T$ (see Figure 1).
The lengths of the graph are given by the
heights of the corresponding cylinders. It then follows that there
is a continuous energy-minimizing map $f:\SR\to T$ in the homotopy
class of the projection along leaves $\pi:\SR\to T$: this follows
because $T$ is non-positively curved. (In this form, this result is
a special case of the far deeper result of Schoen in \cite{GS}. In the
present situation, a more elementary proof is possible, and is
effectively in the literature: we shall include a proof in the
appendix for the
sake of completeness.) The Hopf differential $\Phi$ is then
holomorphic. Moreover, if the vector $V$ is tangent to a regular
vertical trajectory of $\Phi$, then $Vf=0$, from which it follows
that $f$ is constant along leaves of the vertical foliation of
$\Phi$. Since the homotopy class of $f$ is non-constant, the leaves
must be nowhere dense, and hence closed. It also follows that the
map $f$ is given precisely by projecting along the vertical leaves of
the foliation of $\Phi$ to the (vertical) leaf space of $\Phi$.
Thus both $f$ and $\pi$ are defined by projecting along leaves of a
(singular) foliation of $\SR$ by closed leaves: we show that these
two foliations are Whitehead equivalent (defined below) by using
that $f$ and $\pi$ are homotopic to show that the vertical foliation
of $\Phi$ has the same $k$ cylinders of the same heights as $\SF$.
This concludes the argument in the case where the leaf space $T$ of
$\SF$ is a boundaryless graph. In the case where $T$ has boundary,
we approximate $\SF$ by a foliation $\SF_\e$ whose graph $T_\e$ is
boundaryless ($T_\e$ is constructed from $T$ by attaching a small
loop at each boundary point of $T$)
and obtain the quadratic differential representative
of $\SF$ as a limit.

We organize this paper as follows. In the second section, we define
our terms, set our notation and recall some background information.
In \S3, we prove the main result. The paper concludes with an
appendix giving a reasonably elementary proof of the existence of
harmonic maps from surfaces to boundaryless graphs.

\vskip1cm
\sub{\S2 Notation and Background}

\bpage
\sub{2.1 Quadratic Differentials} A measured foliation $(\SF,\mu)$ on
a differentiable surface $F^2$ consists of a foliation of $F$ with
isolated singularities so that the foliation at the singularities
has $k$-pronged singularities, and a translation-invariant measure
$\mu$ supported on arcs transverse to the foliation $\SF$. We will
be interested in foliations all of whose leaves are closed. We will
say that two measured foliations are Whitehead equivalent if the
foliations can be forced to agree by a series of Whitehead moves and
measure preserving isotopies of the underlying surface. (See
\cite{FLP}).

A holomorphic quadratic differential $\Phi$ on a Riemann surface
$\SR$ is a tensor given locally by an expression $\Phi=q(z)dz^2$,
where $z$ is a conformal coordinate on $\SR$ and $q(z)$ is
holomorphic. Such a quadratic differential $\Phi$ defines a
measured foliation in the following way. The zeros $\Phi^{-1}(0)$
of $\Phi$ are well-defined; away from these zeros, we can choose a
canonical conformal coordinate $\z=\int^z\sqrt\Phi$ so that
$\Phi=d\z^2$. The local measured foliations
$(\{\Ree\z=\const\}, d(\Ree\z))$ then piece together to form a measured
foliation known as the vertical measured foliation of $\Phi$. There
is a corresponding horizontal measured foliation constructed out of the
local foliations $(\{\im\z=\const\}, d\im\z)$.

Relative to the vertical foliation of a holomorphic quadratic
differential, the Riemann surface decomposes into two types of
dense open domains: cylindrical domains, foliated by closed curves
all freely homotopic to the same closed curve, and spiral domains,
in which all leaves are non-compact and dense (see \cite{Str3} and
\cite{Gar} for details). We will be interested in quadratic
differentials, and more generally, measured foliations, whose foliations
consist
entirely of cylindrical domains: these are uniquely determined (up
to equivalence, in the measured foliation case) by the free
homotopy classes of the core curves and the heights of the
cylinders (\cite{Str3; \S20}).

\sub{2.2 Energy-Minimizing Maps} Given a map $w: \SR\to(T,h)$ from a
Riemann surface $\SR$ to a locally finite Riemannian complex $T$, we
define the energy form to be the tensor
$edz\ox d\bar z=(\|w_*\p_z\|_h^2+\|w_*\p_{\bar z}\|_h^2)dz\ox
d\bar z$; the energy of the map is $E=\int edz\land d\bar z$. An
energy minimizing map is a minimum for this functional in a
homotopy class. We define the Hopf differential $\Phi$ for a map
$w: \SR\to T$ by
$\Phi=\Phi dz^2=4\left<w_*\p_z, w_*\p_z\right>_hdz^2$.
Note that $\|\Phi\|=\|\Phi\|_{L^1}<2E$

R\.~Schoen
has emphasized \cite{S} that a map for which the energy functional
is stationary under reparametrizations of the domain has a Hopf
differential which is holomorphic: one uses suitable domain
reparametrizations to show that the Hopf differential satisfies the
Cauchy-Riemann equations weakly, and then Weyl's lemma forces the
Hopf differential to be (strongly) holomorphic. We observe that in
this argument, the range manifold may be singular.

The vertical and horizontal foliations of the Hopf differential for
$w: \SR\to T$ integrate the directions of minimal and maximal
stretch of the gradient map $dw$, for smooth energy minimizing maps
$w: \SR\to T$.

\vskip1cm
\sub{\S3 Main Result}

\bpage
Let $\SR$ be a Riemann surface and choose $k$ mutually disjoint
homotopically non-trivial Jordan curves $\vecs\g k$ on $\SR$ no
pair of which are freely homotopic (an admissable system; see
\cite{Str3; \S2.6}). Also choose $k$ positive numbers $\vecs hk$ to
serve as heights, and construct the measured foliation $(\SF,\mu)$
on $\SR$ consisting entirely of cylinders $C_j$ of height $h_j$
with core curves $\g_j$ and compact singular curves. The leaf space
$T$ of $\SF$ is a compact $1$-complex, possibly with boundary
and/or a finite number of finite valence vertices corresponding to
the singular curves (see Figures 1 and 2);
$T$ inherits a metric $d=d_T=\pi_*\mu$ by
pushing forward the measure $\mu$ by the natural projection
$\pi: \SR\to T$ along the leaves. We reprove

\proclaim{Theorem 1} {\rm (\cite{Ren}, \cite{HM})} There is a
holomorphic quadratic differential $\Phi$ on $\SR$ whose vertical
foliation is Whitehead equivalent to $(\SF,\mu)$.
\endproclaim

\sub{Remark} The novelty of our proof is that $\Phi$ will emerge as
the Hopf differential of the (unique) energy minimizing map
$f: \SR\to(T,d)$.

\sub{Proof of Theorem 1} We first consider the case when
$\p T=\emptyset$. Then by Theorem~2, there is a continuous energy
minimizing map $f=f_\pi: \SR\to T$ homotopic to $\pi: \SR\to T$.
The Hopf differential $\Phi=(f^*d)^{2,0}$ is then holomorphic, and
does not vanish everywhere, because this would imply that locally $f$
was either a continuous conformal map between a Riemann surface and
a graph, which is absurd, or constant, which is excluded by the
homotopy between $f$ and $\pi$. We consider its vertical measured
foliation. Away from $\Phi^{-1}(0)$, we consider the natural
coordinate $\z=\int^p\sqrt\Phi$ and we write $\z=\xi+iy$. Of
course, in these coordinates we have $\Phi=1d\z^2$ so that
from our definition of the Hopf differential, we have the two real
equations, $\|f_*\p_\xi\|_d^2-\|f_*\p_\eta\|^2=1$ and
$\left<f_*\p_\xi, f_*\p_\eta\right>_d=0$, at least away from those
neighborhoods $\SU$ in $\SR$ which are mapped onto three or higher
pronged neighborhoods of the singularities of $T$. Since away from
the singularities of $T$, the graph $T$ is a $1$-manifold, we can compare
those two equations to find,
on those neighborhoods $\SU\subset\SR$, that
$f_*\p_\xi\|f_*\p_\eta$. We conclude that $f_*\p_\eta=0$ and
$f_*$ maps the vector $\p_\xi$ onto a unit vector tangent to $T$;
note from the uniform continuity of $f$, that in fact this argument
applies to neighborhoods of all but the finitely many arcs in $\SR$
representing preimages of a singularity of $T$. Thus, the map $f$
is in fact a projection along the leaves of the vertical foliation
of $\Phi$. Now, a leaf of a vertical foliation of a quadratic
differential is either compact and contains a point of
$\Phi^{-1}(0)$, or is closed and avoids $\Phi^{-1}(0)$, or is dense
in some open set in $\SR$
(see \cite{Str3}). If the latter occured, this would force $f$ to
be constant on some open set in $\SR$; this again forces the
holomorphic tensor $\Phi$ to vanish on $\SR$, a contradiction. The
classification of trajectories of a holomorphic quadratic
differential on $\SR$ then implies that the vertical foliation of
$\Phi$ is composed of $\ell$ cylindrical domains, bounded by some
compact trajectories emanating from the (finite) set $\Phi^{-1}(0)$.

Recall the set $\SS$ of simple closed curves on $\SR$. Now
$f: \SR\to T$ is homotopic to $\pi: \SR\to T$, so the set
$A_f=\{[\g]\in\SS\mid\inf\lm_{\g\in[\g]}\ell_{d_T}(f(\g))=0\}$,
of curves with representatives whose image under $f$ can have
arbitrarily small $d_T$-length, agrees with $A_\pi$. Within
$A_f=A_\pi$, we can distinguish
$B_f=B_\pi=\{[\g]\in A_f\mid i([\g],[\b])=0$ $\forall[\b]\in A_f\}$.

\proclaim{Lemma} $B_\pi=\{[\g_1],\dots,[\g_k]\}$, the set of core
curves of the foliation $\SF$.
\endproclaim

\sub{Proof of Lemma} It is clear that
$\{[\g_1],\dots,[\g_k]\}\subset B_\pi$, since each $\g_j$ is
clearly in $A_\pi$ by construction, and any element $[\b]\in\SS$
having an essential intersection with $[\g_j]$ would be mapped by
$\pi$ onto a portion of $T$ which would necessarily include the
segment of $d_T$-length $h_j$ corresponding to the projection of the
cylinder $C_j$.

On the other hand, the complement of $\{[\g_1],\dots,[\g_k]\}$ in
$A_\pi$ consists of classes of curves representable by curves in a
singular trajectory of $\SF$. Consider such a curve class $[\b]$
and consider a representative $\b$ of $[\b]$ in a small
neighborhood $N$ of the singular trajectory. Since $[\b]\neq[\g_j]$
for any $j$, we may take $[\b]$ to be non-boundary parallel. But
then there exists a class $[\a]$ of simple closed curves in $N$
with $i([\a], [\b])\neq0$. (This follows, once we consider a
representative $\b$ of $[\b]$ in $N$: the open manifold $N$ is
planar, so $\b$ must separate $N$. But,
since $\b$ is not boundary parallel, neither component of the complement
of $\b$ in $N$ is a
cylinder, so we can find simple arcs in each component which are
homotopically non-trivial, rel~$\b$. We then connect these arcs to
get the desired curve.) Naturally, then, the class $[\a]$ has a
representative $\a$ in the singular trajectory, so that
$\ell_{d_T}(f(\a))=0$, and so $[\a]$, thought of as an element of
$\SS$, has $[\a]\in A_\pi$. We conclude that $[\b]\in A_\pi-B_\pi$,
proving the lemma. \qed

\sub{Conclusion of the proof of Theorem 1} The lemma, applied to
both $f$ and $\pi$, shows that the vertical foliation for $\Phi$
and $\SF$ share the same core cylinders. The heights $\{h_j\}$ of
the cylinders can be determined by computing
$\inf\lm_{\a\in[\a]}\ell_{d_T}(f(\a))=\inf\lm_{\a\in[\a]}\ell_{d_T}
(\pi(\a))=L([\a])$ for a sufficient number of curve classes
$[\a]\in\SS$, and solving for heights using the identity
$L[\a]=\sum\lm_j i([\a],[\g_j])h_j$.

Since the vertical foliation of $\Phi$ and the foliation $\SF$ have
the same core curves and heights, they are Whitehead equivalent as
desired.

In the case where the graph $T$ has boundary $\{p_1,\dots,p_\ell\}$
we proceed as follows. The pre-image of a boundary point consists
of a singular trajectory. We consider, for each such boundary point
$p_j$, a simple closed curve $\a_j$ contained in the corresponding
singular trajectory (the trajectories are required to be
homotopically non-trivial), and a vector $\vec\e$ of small numbers
$\vec\e=(\e_1,...,\e_l)$. Our
construction gives that the set of curves
$\G=\{[\g_1],\dots,[\g_k]$, $[\a_1],\dots,[\a_\ell]\}$ have
mutually disjoint representatives, no pair of which are freely
homotopic. We construct the foliation $\SF_{\vec\e}$ with core curves
from $\G$ and corresponding heights $\{h_1,\dots,h_k$,
$\e_1,\dots,\e_\ell\}$; we have chosen the curves $\a_j$ so that
the corresponding cylinders $C'_j$ will have both boundary
components on the boundary of a cylinder $C$ corresponding to one
of the $\g_i$, in fact to the $\g_i$ whose cylinder projected under
$\pi$ to a neighborhood of the boundary point $p_j$. We see that
the graph of the leaf space is boundaryless, so that our previous
argument yields a quadratic differential $\Phi_{\vec\e}$ with core
curves $\{[\g_1],\dots,[\g_k]$, $[\a_1],\dots,[\a_\ell]\}$ and
heights $\{h_1,\dots,h_k$, $\e_1,\dots,\e_\ell\}$.

We claim that $\|\Phi_{\vec\e}\|$ is uniformly bounded: then in
that case, as we let
$\vec\e=(\e_1,\dots,\e_\ell)\to0$, the differentials $\Phi_{\vec\e}$
converge uniformly to a quadratic differential with the prescribed
core curves and heights.

To see that $\|\Phi_{\vec\e}\|$ is uniformly bounded, we return to the
problem of finding an energy minimizing map
$w_{\vec\e}: \SR\to(T_{\vec\e},d_{\vec\e})$ where $(T_{\vec\e},d_{\vec\e})$
is the metric graph obtained from the leaf space of the measured
foliation $(\SF_{\vec\e},\mu_{\vec\e})$.  We construct a competitor
$w:\SR\to(T_{\vec\e}, d_{\vec\e})$ to the energy minimizing map
$w_{\vec\e}$ by fixing a vector $\vec\k=(\k_1,...\k_l)$ first and
considering the energy minimizing map
$w_{\vec\k}: \SR\to(T_{\vec\k},d_{\vec\k})$. Now, we have constructed
$T_{\vec\k}$ and $T_{\vec\e}$ to be homeomorphic, with a canonical map
between the vertices of the graph.  We construct an affine map
$A_{\vec\k,\vec\e}: (T_{\vec\k}, d_{\vec\k)}\to(T_{\vec\e},d_{\vec\e})$
given by extending that map of vertices of the graph to the
one-complex in an affine way. The composed map
$A_{\vec\k,\vec\e} \circ w_{\vec\e}: \SR\to(T_{\vec\e},d_{\vec\e})$ has
total energy $E_{\vec\k}$ exceeding $E(w_{\vec\e})$ because $w_{\vec\e}$
is energy minimizing; moreover, since for $\vec\e$ with $\e_i<\k_i$ for
$i=1,...l$, the map $A_{\vec\k,\vec\e}$ differs
from an isometry only by being a contraction on the branches
of $T_{\vec\k}$ dual to the
core curves $\{[\a_1],...,[\a_l]\}$, we see that
we have the uniform bound $E_{\vec\k}<E_0$.
This yields the estimate
$\|\Phi_{\vec\e}\|<2E(w_{\vec\e})<2E_{\vec\k}<2E_0$ (where we
recall the first inequality from \S2), proving the claim.

\vskip1cm
\sub{\S4 Appendix}

\bpage
We prove, with an argument (see \cite{J1; \S4.1}, and \cite{GS})
that is already almost entirely in the literature,

\proclaim{Theorem 2} Let $\phi: \SR\to(T,d)$ be a map from a
Riemann surface $\SR$ to a graph $T$, where $T$ is equipped with a
metric $d$; suppose $\phi\in C^0\cap H^1(\SR,T)$. Then there is a
continuous energy minimizing map $u: \SR\to(T,d)$ homotopic to
$\phi$ with the modulus of continuity of $u$ estimable in terms of
$E(\phi)$ and the modulus of continuity of $\phi$.
\endproclaim

\sub{Remark} We include this proof for the purpose of displaying
the analytical underpinnings of our approach to Strebel's theorem,
not for any claim of novelty.

\sub{Proof} The plan is to prove the result first locally, and then
to piece together the local harmonic maps into a well controlled
sequence of maps which tend towards a minimizer $u: \SR\to(T,d)$.

\sub{4.1} So consider first a continuous map of the circle
$g: \p\Delta\to(T^*,d)$ where $g$ admits an extension
$\bar g: \p\Delta\to(T^*,d)$ of finite energy and the image of $g$
is a simply connected closed subgraph $T^*\subset T$, for instance
the intersection of $T$ with a ball,
$T^*=T \cap B(p,r)$. We claim that there exists a harmonic map
$h: \Delta\to(T,d)$ with boundary values $g$, and that $h$
minimizes the energy with respect to these boundary values. We also
claim that the modulus of continuity of $h$ can be estimated in
terms of $E(\bar g)$ and the modulus of continuity of $g$. To see
this, take a minimizing sequence $\{v_i\}$ for the energy in
$V=\{v\in H^1(\Delta, B(p,r^*))$, $v\bigm|_{\p\Delta}=g\}$ where
$r<r^*$. (The point here is that our space $V$ makes perfectly good
sense for a singular target being a graph, since we can embed the
simply connected graph $T^*$ isometrically in $R^2$, defining
$H^1(\Delta, B(p,r^*))=\{v\in H^1(\Delta,\BR^2)$;
$v(z)\in B(p,r^*)\subset T^*\subset\BR^2$ a.e. $z\in\Delta\}$ and
$v=\hat v$ on $\p\Delta$ if $v-\hat v\in H^1_0(\Delta,\BR^2)$
where $H^1_0(\Delta,\BR^2)$ is the $H^1$-norm closure of smooth,
compactly supported $\BR^2$-valued maps on $\Delta$.) Now, such a
minimizing sequence $\{v_i\}$ has a subsequence converging weakly
in $H^1$, and the limit $h$ minimizes energy in its class because
of the lower semicontinuity of the energy functional. Since the set
$\{v\in H^1(\Delta,\BR^2)$, $v\bigm|_{\p\Delta}=g\bigm|_{\p\Delta}$
is a closed affine subspace of $H^1(\Delta,\BR^2)$, it is weakly
closed, and we conclude that
$h\bigm|_{\p\Delta}=g\bigm|_{\p\Delta}$. We can estimate the
modulus of continuity of $h$ as follows: the Courant-Lebesgue lemma
(\cite{J1; Lemma~3.1.1}) provides that for each $x\in\SR$ and $\e>0$
there is a $\rho$ depending only upon $\e$, $E(\bar g)$ and the
modulus of continuity of $g$ and $q\in T^*$ so that
$h(\p B(x,\rho))\subset B(q,\e)\cap T^*$. But then, since
$B(q,\e)\cap T^*$ is convex, the maximum principle forces
$h(B(x,\rho))\subset B(q,\e)\cap T^*$. Thus we have shown the
continuity of $h$ and estimated its modulus of continuity.

\sub{Remark} Of course, when $T^*\subset T$ is an interval, $h$ is
given classically by the Poisson integral formula, so here we are
really only interested in $T^*$ being a non-trivial graph, for
instance a ``$Y$''. It would be interesting to have a more explicit
description of the map $h: \ov\Delta\to T^*$ in this case, in terms
of the boundary values $h\bigm|_{\p\Delta}: \p\Delta\to T^*$.

\sub{4.2} To complete the argument, we choose a $\d_0$ sufficiently
small so that the Courant-Lebesgue lemma forces a harmonic map $\phi$ with
energy $E(\phi)$ of the ball of radius $\d<\d_0$ to be contained in
a simply connected subgraph $T^*\subset T$; we want to be able to
apply the construction of the previous paragraph 4.1.

Choose $\d<\d_0$, and cover $\SR$ by a finite number $M$ of balls
$B(x_i,\d/2)$, $i=1,\dots,M$; here we have chosen some suitable
background metric on $\SR$. Following \cite{J1; Theorem~4.1.1}, we
let $u_n$ be a continuous energy minimizing sequence of maps
homotopic to $\phi$; we may as well assume that $E(u_n)<E(\phi)$ so
that we have control on $u_n(B(x_i,\d/2))$. Thus, by the
Courant-Lebesgue lemma and our choice of constants, for every $n$,
we can find $r_{n,1}$, with $\d<r_{n,1}<\d^{1/2}$ and $p_{n,1}\in N$
so that $u_n(\p B(x_1,r_{n,1}))\subset B(p_{n,1},r^*)$ for $r^*$ so
small that $B(p,r^*) \subset T$ is always simply connected. Thus, we can
invoke our construction for harmonic maps of balls into simply
connected graphs: we replace $u_n\bigm|_{B(x_1,r_{n,1})}$ by the
harmonic map $h_{x_1,n}: B(x_1,r_{n,1})\to B(p_{n,1},r)$ whose
boundary values $h_{x_1,n}\bigm|_{\p B(x_1,r_{n,1})}$ agree
pointwise with those $u_n\bigm|_{\p B(x_1,r_{n,1})}$ of $u_n$.

We can assume that $r_{n,1}\to r_1$ and using the estimates on the
modulus of continuity of the harmonic maps $h_{x_1,n}$ of simply
connected domains, we can take the replaced maps, say $u_{n,1}$, to
converge uniformly on $B(x_1,\d-\eta)$ for any $0<\eta<\d$.
Naturally $E(u_{n,1})<E(u_n)$.

We repeat the argument to find $r_{n,2}$ with $\d<r_{n,2}<\sqrt\d$
and $u_{n,1}(\p B(x_2,r_{n,2}))\subset B(p_{n,2},r)$ and replace
$u_{n,1}$ on $B(x_2,r_{n,2})$ by the appropriate solution to the
Dirichlet problem; we denote the new replaced maps by $u_{n,2}$,
and again assume that $r_{n,2}\to r_2$.

Now, in the first replacement step, $u_{n,1}$ converged uniformly on
$B(x_2,r_2)\cap B(x_1,\d-\eta/2)$, and thus the boundary values for
our second replacement step converge uniformly on
$\p B(x_2,r_{n,2})\cap B(x_1,\d-\eta/2)$. Of course, we have a uniform
estimate for the modulus of continuity between small disks and
portions of graphs (even after allowing for the difference between
our background metric and the Euclidean metric on the disk): this
means that we can assume that the maps $u_{n,2}$ converge uniformly
on $B(x_1,\d-\eta)\cup B(x_2,\d-\eta)$ for $0<\eta<\d$. Of course,
the act of replacing lowers energy, and we conclude that
$E(u_{n,2})\le E(u_{n,1})\le E(u_n)$. We repeat the replacement
argument in this way obtaining a family of maps $u_{n,M}$ with
$E(u_{n,M})\le E(u_n)$, and which converge uniformly on all balls
$B(x_i,\d/2)$, hence on $\SR$ since
$\SR\subset\bigcup\lm^M_{i=1}B(x_i,\d/2)$.

We let $u$ denote the limit of $u_{n,M}$ as $n\to\infty$ (recall that
$M$ is a
fixed number depending only upon $\d$). Note that since replacement
on disks will not affect the homotopy class of the maps $u_{n,i}$,
the uniform convergence of $u_{n,M}$ to $u$ forces $u$ to be not
only continuous but also
homotopic to the given $\phi: \SR\to T$. Naturally also, since $u_n$
is a minimizing sequence, so is $u_{n,M}$; we note that since
$E(u_{n,M})\le E(\phi)$, the maps $u_{n,M}$ converge weakly in
$H^1(\SR,T)$ (see above comments on proper definitions), and by the
lower semicontinuity of the energy functional with respect to weak
$H^1$ convergence, we get that the limit $u$ of the energy
minimizing sequence $u^M_n$ minimizes energy within the homotopy
class of $\phi: \SR\to T$.

\sub{Remarks} 1) When the image of $u$ is contained in an embedded
interval in the graph, it follows immediately that $u$ is a
harmonic function, and hence real analytic. 2) It is possible to
start with $\phi: \SR\to T$ being a projection as in the application
(Theorem~1), and then argue that all replacements involve only
Whitehead moves to the resulting Hopf differential vertical
foliation. This yields the desired Jenkins-Strebel differential
possibly more constructively.

\end
\vskip1cm
\Refs

\widestnumber\key{Str3}

\ref
\key FLP  \by A. Fathi, F. Laudenbach, and V. Poenaru
\paper Traveaux de Thurston sur les Surfaces
\jour Asterisque     \yr1979     \pages66--67
\endref

\ref
\key Gar  \by F.P. Gardiner
\book Teichm\"uller Theory and Quadratic Differentials
\publ Wiley   \publaddr New York   \yr1987
\endref

\ref
\key GS    \by  M. Gromov and R. Schoen
\paper Harmonic Maps into Singular Spaces and $p$-adic Superrigidity
for Lattices in Groups of Rank One
\paperinfo to appear in Publ. IHES
\endref

\ref
\key HM    \by  J. Hubbard and H. Masur
\paper Quadratic Differentials and Foliations
\jour Acta Math.   \vol142     \pages221--224 (1979)
\endref

\ref
\key  Je   \by J.A. Jenkins
\paper On the Existence of Certain General Extremal Metrics
\jour Ann. of Math.  \vol66   \pages440--453   \yr1957
\endref

\ref
\key  J1   \by J. Jost
\book Two Dimensional Geometric Variational Problems
\publ Wiley  \publaddr West Sussex, England      \yr1991
\endref

\ref
\key  J2   \bysame
\paper Existence proofs for harmonic mappings with the help of a
maximum principle
\jour Math. Z.  \vol184   \yr1983  \pages489--496
\endref

\ref
\key  Mo   \by  C.B. Morrey
\paper The Problem of Plateau on a Riemannian Manifold
\jour Ann. of Math.  \vol49   \yr1978  \pages807--851
\endref

\ref
\key  Ren   \by  H. Renelt
\paper Konstruktion gewisser quadratischer Differentiale mit
Hilfe von Dirichletintegralen
\jour Math. Nachr.  \vol73   \yr1976  \pages125--142
\endref

\ref
\key S  \by R. Schoen
\book Analytic Aspects of the Harmonic Map Problem
\bookinfo in Seminar on Nonlinear Partial Differential Equations,
S.S. Chern, ed., MSRI Publ 2, Springer, New York, 321--358
\endref

\ref
\key Str1   \by  K. Strebel
\book \"Uber quadratische Differentiale mit Geschlossenen Trajektorien
und extremale quasikonforme Abbildungen
\bookinfo Festband zum 70. Geburtstag von Rolf Nevanlinna (1965/1966)
Springer Verlag, 105--127
\endref

\ref
\key Str2   \bysame
\paper Bemurkungen \"uber quadratische Differentiale mit
geschlossenen Trajektorien
\jour Ann. Acad. Sci. Fenn. A.I.  \vol405  \yr1967  \pages1--12
\endref

\ref
\key Str3   \bysame
\book Quadratic Differentials
\publ Springer   \publaddr Berlin  \yr1984
\endref

\endRefs
